\documentclass[twoside,12pt]{article}
\usepackage{amsmath,amsthm,amssymb,amscd}
\usepackage[matrix,arrow,curve]{xy}
\title {Pseudo-isotopy classes of diffeomorphisms of the unknotted pairs
$(S^{n+2},~S^n)$ and $(S^{2p+2},~S^p\times S^p)$ }

\author{Nikolai A. Krylov
\\
\small Siena College, School of Science\\
\small 515 Loudon Road\\
\small Loudonville, New York 12211\\
\small e-mail:~nkrylov@siena.edu}

\date {}

\newtheorem{thm}{Theorem}
\newtheorem{cor}{Corollary}
\newtheorem{lem}{Lemma}

\newtheorem{prop}{Proposition}

\def\int            {\mathbb Z}

\def\rea            {\mathbb R}
\def\com            {\mathbb C}

\def\F              {{\cal F}}

\def\mcgM           {\pi_0 {\rm Diff}(M) }

\def\Gv2            {\Gamma_V(2)}

\def\lra            {\longrightarrow}

\def\hra            {\hookrightarrow}
\def\lmt            {\longmapsto}

\def\vf             {\varphi}

\begin{document}

\maketitle

\parskip=5mm

\begin{abstract}
{We consider two pairs: the standard unknotted $n$-sphere in
$S^{n+2}$, and the product of two $p$-spheres trivially embedded in
$S^{2p+2}$, and study orientation preserving diffeomorphisms of
these pairs. Pseudo-isotopy classes of such diffeomorphisms form
subgroups of the mapping class groups of $S^n$ and $S^p\times S^p$
respectively and we determine the algebraic structure of such
subgroups when $n>4$ and $p>1$.}
\end{abstract}

\noindent {\bf Keywords}: Isotopy of diffeomorphisms;
quadratic forms; non-spherical knots\\
{\bf 2000 Mathematics Subject Classification}: 57N37; 57Q45; 57R50

\section{Introduction}

Let $M^n$ be a closed, oriented, smooth manifold of dimension $n$.
We denote the group of all orientation preserving diffeomorphisms of
$M$ by ${\rm Diff}(M)$. Recall that two elements of this group, say
$f$ and $h$, are called pseudo-isotopic if there exists an
orientation preserving diffeomorphism $F\in {\rm Diff}(M\times
[0,1])$ such that $F(x,0)=f(x)$ and $F(x,1)=h(x)$ for all $x\in M$.
The set of pseudo-isotopy classes forms a group which we will call
{\sl the mapping class group of $M$} and denote by $\mcgM$. For some
special types of manifolds this group has been computed by several
authors, see for example \cite{Fried}, \cite{Kreck}, \cite{Krylov1},
\cite{Levine2}, \cite{Sato}, \cite{Turner}, \cite{WellsRob}. On the
contrary, the classical mapping class group ${\cal M}_g$ that is the
group of the isotopy classes of orientation preserving
diffeomorphisms of an orientable surface $F^2$ of genus $g$ has been
studied extensively and there are also results on the so-called {\sl
spin mapping class group} (\cite{Harer1}, \cite{Hirose},
\cite{Montes}). This last one is a subgroup of ${\cal M}_g$ that
consists of diffeomorphisms preserving the Rokhlin quadratic form
$\rho: H_1(F^2;\int_2)\lra \int_2$ (see next section for the
definition). This spin mapping class group has a geometrical
meaning: it consists of the diffeomorphisms that extend to
diffeomorphisms of the 4-sphere, if one assumes that the closed
2-surface is trivially embedded in $S^4$ (see \cite{Hirose} and
\cite{Montes}). From this point of view, there are interesting
results of S. Hirose (Theorem A. in \cite{Hirose2}) and Z. Iwase
(Proposition 3.1 in \cite{Iwase}) who have shown that for a
non-trivially embedded 2-torus $T^2\hra S^4$ the isotopy classes of
diffeomorphisms of $T^2$ which extend over $S^4$ form a proper
subgroup of the spin mapping class group of the torus. One way to
obtain such a non-trivial embedding is to take a classical torus
knot $(S^3,~K)$, choose a 3-disk $D^3\hra S^3\setminus (K\times
D^2)$ and construct the {\it spun $T^2$-knot} in $S^4$ using an
analog of the spinning construction:
$$
\left(S^4,~S(K)\right):=\left((S^3,~K)\setminus {\rm
int}D^3\right)\times S^1 \bigcup_{\rm id} S^2\times D^2
$$

One could consider $2-$surfaces not only in $S^4$ but in any other
closed orientable $4-$manifold and ask similar questions. An
attractive family of $2-$surfaces in a $4-$manifold is the family of
non-singular algebraic curves, say in $\com P^2$. Theorem 4.3 of
\cite{Hirose3} claims that for the smooth algebraic curves of odd
degree greater than four, the isotopy classes of diffeomorphisms
extendable over $\com P^2$ form a proper subgroup of the
corresponding mapping class group.

Here we are interested in a higher dimensional analog. Let
$M^n,~n\geq 4$ be a closed oriented locally unknotted submanifold of
the standard sphere $S^{n+2}$, in other words, consider an
$n-$dimensional locally flat differentiable $M^n-$knot
$(S^{n+2},~M^n)$. The following questions naturally rise:
\begin{itemize}
\item Given such a knot, which elements of the mapping class
group $\pi_0{\rm Diff}(M^n)$ have representatives that extend to
orientation-preserving diffeomorphisms of the ambient sphere
$S^{n+2}$?
\item How does the set of such elements depend on the
knot type?
\end{itemize}
It should be mentioned that similar questions appear in the study of
the homotopy type of the complement to an algebraic hypersurface in
the complex space. For the details, the reader is referred to a
recent work of A. Libgober (see \cite{Libgober}, \S3.2), and to \S 5
of \cite{Dolgachev}.

In this paper we answer the first question for two ``trivial" knots:
$(S^{n+2},~S^n)$ with $n\geq 5$, and $(S^{2p+2},~S^p\times S^p)$
with $p\geq 2$.

Let us define first the subgroup of the mapping class group of $M$,
which we will be dealing with. Suppose $M^n,~n\geq 5$ is a closed
simply-connected manifold embedded in the standard sphere
$S^{n+m},~m\geq 1$. Consider two orientation preserving
diffeomorphisms of $M$, say $\phi_1$ and $\phi_2$, which are
isotopic. If we assume that $\phi_1$ extends to an orientation
preserving diffeomorphism of the ambient sphere, then the isotopy
extension theorem (see, for example, \cite{Hirsch}, ch. 8, Theorem
1.5) guarantees that $\phi_2$ can also be extended to a
diffeomorphism of $S^{n+m}$. Moreover, it follows from the disk
theorem (\cite{Hirsch}, ch. 8, Theorem 3.1) that we can always
assume the extension is isotopic to the identity map of the ambient
sphere (by changing the extension on a small disk embedded far away
from $M$). This implies that for a $M^n$-knot $(S^{n+m},~M)$ the
following subgroup of $\mcgM$, denoted by ${\cal E}(S^{n+m},~M)$, is
well defined (cf. definition of ${\cal E}(S^4,~K)$ in
\cite{Hirose}).
$$
{\cal E}(S^{n+m},~M):=\left\{[\phi]\in\mcgM ~|~{\rm there~ exists~}
\Phi\in {\rm Diff}(S^{n+m})\\ {\rm ~s.~ t.~} \Phi|_M=\phi \right\}
$$
When $n\geq 5$, we can replace an isotopy by pseudo-isotopy and the
freedom to make such a change is guaranteed by the result of J. Cerf
\cite{Cerf}.

In the next section we will use a higher dimensional analog of the
Rokhlin quadratic form introduced by J. Levine in \cite{Levine1} to
generalize the result of J. Montesinos (see \cite{Montes}) who has
shown that for the trivially embedded torus $T^2\hra S^4$ the group
${\cal E}(S^4,~T^2)$ consists of the matrices $
\begin{pmatrix}
a & b\\
c & d\\
\end{pmatrix}\in {\rm SL}(2,\int)$ where both products $a\cdot b$
and $c\cdot d$ are even integers. In the last paragraph we study
diffeomorphisms of $S^p\times S^p$ which act trivially on the
homology and extend over the ambient sphere $S^{2p+2}$. Isotopy
classes of such diffeomorphisms form a subgroup of ${\cal
E}(S^{2p+2},~S^p\times S^p)$, which will be denoted by ${\cal
SE}(S^{2p+2},~S^p\times S^p)$. We will write $\phi\sim \vf$, and
$G\cong H$ to say that diffeomorphisms $\phi$ and $\vf$ are (pseudo)
isotopic, and groups $G$ and $H$ are isomorphic respectively. All
homology and cohomology groups will have integer coefficients unless
otherwise stated.

\section{Levine - Rokhlin quadratic form.}

In this part we discuss a higher dimensional analog of the Rokhlin
quadratic form and show that if a diffeomorphism $f\in {\rm
Diff}(M^n)$ extends to a diffeomorphism of the ambient sphere
$S^{n+2}$, then $f_*$ (the induced homology map) commutes with
this quadratic form $q$. Definition of such a form $q$ is based on
the following construction of Levine (\cite{Levine1}, \S3).

Let us denote the connected sum of $k$ copies of $S^p\times S^p$ by
$M^n$ ($p\geq 3$ is to be odd in this section) and assume that $M^n$
is somehow embedded in $S^{n+2}$. Notice that since $M$ is at least
$2$-connected, the normal bundle
$\nu(S^{n+2},M^n)\simeq\varepsilon^2$ is trivial for any embedding
$M^n\hra S^{n+2}$. Hence we can choose a framing $\F=(f_1,f_2)-$two
orthonormal vector fields on $M$ to get a framed submanifold
$(M,\F)$ of $S^{n+2}$. Moreover, since $M$ is simply connected the
set of homotopy classes $[M,~SO(2)]$ consists of the trivial element
only and therefore the choice of our framing $\F$ is unique up to
homotopy. For such a framed manifold, Levine (see \cite{Levine1},
\S4 and also \cite{Haefliger2}, \S3.2) has defined a function
$$
\phi(M,\F): ~H_p(M;\int)\lra \int_2
$$ which depends only on the
isotopy class of the pair $(M,\F)$ and satisfies the formula
$$
\phi(\alpha + \beta) = \phi(\alpha) + \phi(\beta) + (\alpha \cdot
\beta)_2
$$ where $\alpha,\beta\in H_p(M;\int)$ and
$(\alpha\cdot\beta)_2$ is the mod 2 residue of the intersection
number of $\alpha$ and $\beta$. We now recall this definition of
Levine in a slightly different form which will suit needs of this
paper. Fix an embedding $M^n\hra S^{n+2}$ and a framing
$\F=(f_1,f_2)$ of the normal bundle $\nu(S^{n+2},M^n)$. Take a point
$pt$ outside a tubular neighborhood of $M$ in $S^{n+2}$ and choose a
(positive) framing $\cal G$ of $\tau(S^{n+2})|_{S^{n+2} - pt}$ (the
tangent bundle restricted to the complement of $pt$). Let ${\tt
S}\hra M$ be an embedded $p$-sphere which represents a class $[z]\in
H_p(M;\int)$ and consider the subbundle $E$ of $\tau(S^{n+2})|_{\tt
S}$ which is the Whitney sum of the tangent bundle $\tau({\tt S})$
and the line bundle $f_1|_{\tt S}$. This subbundle has a canonical
framing $\Upsilon_1=(e_1,\ldots,e_p,e_{p+1})$ which comes from the
standard framing of $\rea^{p+1}$ (cf. \cite{Haefliger2}, \S3.2 and
recall that we have fixed an embedding $M^n\hra S^{n+2}$). This
$\Upsilon_1$ together with the field $f_2|_{\tt S}$ gives us a
framing $\Upsilon_2$ of the sum $E\oplus f_2|_{\tt S}$ which is a
$(p+2)$-subbundle of $\tau(S^{n+2})|_{\tt S}$. Comparing this
framing $\Upsilon_2(x)$ with the framing ${\cal G}(x)$ at each point
$x\in {\tt S}$, produces a $(p+2, 2p+2)$-matrix with orthonormal row
vectors, i.e. an element of the Stiefel manifold
$V_{p+2}(\rea^{2p+2})$ of orthonormal $(p+2)$-frames in
$\rea^{2p+2}$. Thus we obtain an element of the homotopy group
$\pi_p(V_{p+2}(\rea^{2p+2}))\cong\int_2$ (for the last isomorphism
see \cite{BrowderB}, Theorem IV.1.12). In our case, this
construction gives a well defined function
$$
\phi(M): ~H_p(M;\int)\lra \int_2
$$ which depends only on the isotopy class of the embedding
$M^n\hra S^{n+2}$.

\begin{prop} {\rm (cf. \cite{Montes}, Proposition 4.1;
\cite{Hirose}, Theorem 1.2)} Let $p=\frac{n}{2}\geq 3$ be odd and
$M^n$ be a closed $(p-1)$-connected manifold embedded into
$S^{n+2}$ and $[f]\in {\cal E}(S^{n+2},~M^n)$. Then $f_*:
H_p(M)\to H_p(M)$ commutes with the function $\phi$.
\end{prop}
\begin{proof}
Let ${\tt S}^p\hra M$ be an embedded sphere representing a cycle
$[{\tt S}^p]\in H_p(M,\int)$. Since diffeomorphism $f$ extends to a
diffeomorphism of the ambient sphere we could compare two
$\Upsilon_2$-framings at each point of $f({\tt S}^p)$. This
comparison would give us an element of $\pi_p(SO(p+2))$. Since the
$(f_1,f_2)-$restrictions of these two framings are homotopic, it
follows from the exact homotopy sequence of the fibration $SO(p)\hra
SO(2p+2)\twoheadrightarrow V_{p+2}(\rea^{2p+2})$ that $\phi([f({\tt
S}^p)]) = \phi([{\tt S}^p])$.
\end{proof}

\noindent \underline{Remark}: If we denote the boundary operator
from this exact homotopy sequence by $\partial$, then
$\partial(\phi([{\tt S}^p]))\in \pi_{p-1}(SO(p))$ is the
obstruction to trivializing the normal bundle of ${\tt S}^p$ in
$M$ (see \cite{Haefliger2}, Lemma 3.4).

Let us now assume that our $2p$-manifold $M^n\simeq (S^p\times
S^p)\# \ldots \# (S^p\times S^p)$ is trivially embedded in
$S^{n+1}$, which is embedded in $S^{n+2}$ as an equator: $M^n\hra
S^{n+1}\hra S^{n+2}$. In particular, we assume that $S^{n+1}$ can be
presented as the union of two copies of the handlebody ${\cal
H}\simeq (D^{p+1}\times S^p)\#\ldots \# (D^{p+1}\times S^p)$ (here
$\#$ stands for the boundary connected sum) along the boundary
$M^n=\partial {\cal H}$:
$$
S^{n+1}={\cal H}_+\bigcup_M-{\cal H}_-.
$$
Taking the tensor product of $H_p(M;\int)$ with $\int_2$ and using
the function $\phi$, we obtain a quadratic form $q:
H_p(M;\int_2)\lra \int_2$ which has Arf invariant zero. For a closed
orientable $4-$manifold $X$ with $H_1(X;\int_2)\simeq 0$ and a
closed orientable embedded 2-surface $F\hra X^4$ which realizes the
element of $H_2(X;\int_2)$ dual to $w_2(X)$, V. Rokhlin ~(see
\cite{Rokhlin}, \S3) has defined a function $\psi: H_1(F;\int_2)\lra
\int_2$ which also satisfies the formula $ \psi(\alpha + \beta) =
\psi(\alpha) + \psi(\beta) + (\alpha \cdot \beta)_2 $. The form
$\psi$ has Arf invariant zero for a trivially embedded surface too.

We now turn to the case where $M\simeq S^p\times S^p$ is standardly
embedded in $S^{2p+2}$. By this we will mean that $M$ in $S^{2p+2}$
is the boundary of a tubular neighborhood of a standardly embedded
$S^p\hra S^{2p+1}$, where $S^{2p+1}$ is an equator of $S^{2p+2}$
(cf. Definition 1.6 of \cite{Lucas}). It then follows from our
Proposition 1 and Propositions 4.2 and 4.3 of \cite{Montes} that if
a diffeomorphism $f\in {\rm Diff}(M)$ extends to a diffeomorphism of
the pair $(S^{n+2},M)$ then the induced automorphism $f_*$ is an
element of the group, which consists of the matrices $
\begin{pmatrix}
d_1 & d_2\\
d_3 & d_4\\
\end{pmatrix}\in {\rm SL}(2,\int)$ where both products $d_1d_2$ and
$d_3d_4$ are even integers. Clearly, any matrix of this type is
congruent modulo 2 either to $ Id=
\begin{pmatrix}
1 & 0\\
0 & 1\\
\end{pmatrix}$ or
$ V :=\begin{pmatrix}
0 & -1\\
1 & 0\\
\end{pmatrix}$ and we will denote such a subgroup of ${\rm SL}(2,\int)$
by $\Gv2 $. One can use the fact that the corresponding projective
group $\Gamma_V(2)/\int_2\cong\int_2 * \int$ is generated by $V$
and $T:=\begin{pmatrix}
1 & 2\\
0 & 1\\
\end{pmatrix}$
(cf. \cite{Wall2}, \S 3) to show that $\Gv2 $ admits the following
presentation $ \Gv2 \cong \langle V,T~|~V^4=Id,~V^2T=TV^2\rangle$.

Let $h:\mcgM \lra {\rm Aut}(H_p(M;\int))$ be the obvious
homomorphism induced by the action on the homology of $M$. We will
denote the restriction of this homomorphism onto the subgroup ${\cal
E}(S^{n+2},~M)$ by $h_{\cal E}$. Evidently, ${\rm Im}(h_{\cal
E})\subset \Gv2 $. Montesinos has shown (see \cite{Montes}, Theorem
5.4) that if $M$ is a standardly embedded torus $S^1\times S^1\hra
S^3\hra S^4$, then ${\cal E}(S^4,~M)\cong \Gv2 $. Here we prove the
following generalization of this result.

\begin{lem} Let $p\geq 3$ be odd and $S^p\times S^p$ be standardly
embedded in $S^{2p+2}$. Then ${\rm Im}(h_{\cal E})= \Gv2 $.
\end{lem}
\begin{proof}
We will show that that there exist diffeomorphisms $\vf_V\in {\rm
Diff}(S^{n+2},~M^n)$ and $\vf_T\in {\rm Diff}(S^{n+2},~M^n)$ such
that $h_{\cal E}[\vf_V]=V$ and $h_{\cal E}[\vf_T]=T$. Let us start
with the matrix $T$ and consider the following commutative diagram
that consists of two exact homotopy sequences of the fibration
$SO(n)\hra SO(n+1) \twoheadrightarrow S^n$
$$
\xymatrix{
& \pi_{p+1}(S^{p+1})\ar[d]_{\partial} \ar[dr]^{\mu=\times 2}\\
\pi_p(SO(p)) \ar[r]_{i_p} & \pi_p(SO(p+1)) \ar@{->>}[d]_{i_{p+1}}
\ar[r]_(.62){j_p} & \pi_p(S^p) \ar[r]_{} & \pi_{p-1}(SO(p)) \\ &
\pi_p(SO(p+2)) }
$$
Since $p$ is odd, map $\mu:=j_p\circ\partial$ is multiplication by 2
(see \cite{BrowderB}, Lemma IV.1.9). It follows from this diagram
that there exists a smooth map $m_T: S^p\lra SO(p+1)$ such that
$j_p([m_T])=\pm2$ and $i_{p+1}[m_T]=0$. Hence there exists a smooth
map (see Smooth Approximation Theorem in \S II.11 of \cite{Bredon})
$\gamma: D^{p+1}\lra SO(p+2)$ that extends the composition of $m_T$
with the inclusion $\iota: SO(p+1)\hra SO(p+2)$, i.e. $\gamma|_{S^p}
= \iota\cdot m_T$. Now we use $\gamma$ together with $m_T$ to define
the following self-diffeomorphism of $S^{2p+2}=(D^{p+1}\times
S^{p+1})\cup (S^p\times D^{p+2})$:
$$
\Phi(x,y):=\left\{
\begin{array}{lcl}
(x,\iota\cdot m_T(x)\circ y)& {\rm if} & (x,y)\in S^p\times D^{p+2}\\
(x,\gamma(x)\circ y)& {\rm if} & (x,y)\in D^{p+1}\times S^{p+1}\\
\end{array}
\right.
$$

Since the restriction of $\Phi$ onto the product $S^p\times S^p\hra
S^p\times S^{p+1}$ (here $x\times S^p$ is the equator of $x\times
S^{p+1}$) is the map $(x,y)\lmt (x,m_T(x)\circ y)$ and
$j_p([m_T])=\pm2$, we can define $\vf_T:=\Phi$ to obtain the
equality $h_{\cal E}[\vf_T]=\pm T$. The product $S^p\times S^p$
bounds a $S^p \times D^{p+1}$ smoothly embedded in $S^{2p+2}$. Then
it follows from Theorem 1.7 of \cite{Lucas} that the pair
$(S^{2p+2},~S^p\times S^p)$ is equivalent to the standard one, where
$S^p\times S^p$ is standardly embedded in the equator $S^{2p+1}$.

For the other case we consider $S^p\times S^p$ standardly embedded
into the unit sphere $S^{2p+2}\subset \rea^{2p+3}$. Using
$\{x_0,x_1,\ldots, x_{2p+2}\}$ as the coordinates in $\rea^{2p+3}$,
we present the equator sphere $S^{2p+1}=\{(x_0,x_1,\ldots,
x_{2p+2})~|~x_0=0~\&~x_1^2+\cdots+x_{2p+2}^2=1\}$ as the union
$S^{2p+1}={\cal H}_1\cup-{\cal H}_2$ where
$$
{\cal H}_1=\{(x_1,x_2,\ldots, x_{2p+2})\in
S^{2p+1}~|~x^2_1+x^2_2+\cdots+x^2_{p+1}\leq \frac{1}{2}\}
$$ and
$$
{\cal H}_2=\{(x_1,x_2,\ldots, x_{2p+2})\in
S^{2p+1}~|~x^2_1+x^2_2+\cdots+x^2_{p+1}\geq \frac{1}{2}\}
$$
and consider a linear map $\Omega: \rea^{2p+3}\lra \rea^{2p+3}$
defined by the square matrix
$$ \Omega :=
\begin{pmatrix}
\begin{tabular}{c|c|c}
$(-1)^p $ & $\begin{matrix} 0 & 0 & \ldots & 0\end{matrix}$ &
$\begin{matrix} 0 & 0 & \ldots & 0\end{matrix}$ \\
\hline
$\begin{matrix} 0\\
0\\
\vdots \\
0\end{matrix}$ & $
\begin{matrix}
0 & 0 & \ldots & 0 \\
0 & 0 & \ldots & 0 \\
\vdots & \vdots & \ddots & \vdots \\
0 & 0 & \ldots & 0
\end{matrix}$ &
$\begin{matrix}
1 & 0 & \ldots & 0 \\
0 & 1 & \ldots & 0 \\
\vdots & \vdots & \ddots & \vdots \\
0 & 0 & \ldots & 1
\end{matrix}$\\
\hline
$\begin{matrix} 0\\
0\\
\vdots \\
0\end{matrix}$ & $\begin{matrix}
-1 & 0 & \ldots & 0 \\
0 & 1 & \ldots & 0 \\
\vdots & \vdots & \ddots & \vdots \\
0 & 0 & \ldots & 1
\end{matrix}$ &
$\begin{matrix}
0 & 0 & \ldots & 0 \\
0 & 0 & \ldots & 0 \\
\vdots & \vdots & \ddots & \vdots \\
0 & 0 & \ldots & 0
\end{matrix}$
\end{tabular}
\end{pmatrix}
$$
of size $(2p+3)\times (2p+3)$ (three blocks in the lowest row have
sizes $(p+1)\times 1$, $(p+1)\times (p+1)$, and $(p+1)\times (p+1)$
correspondingly). $\Omega$ acts on vectors of $\rea^{2p+3}$ by
multiplying a row vector by $\Omega$ on the right.

Since $\Omega$ is an orthogonal matrix which maps the unit disk
$D^{2p+3}$ onto itself preserving orientation
($\rm{det}(\Omega)=(-1)^{2p+2}=1$), restriction of $\Omega$ onto the
unit sphere $S^{2p+2}$ will be an orientation preserving
diffeomorphism as well (see \S4.4 of \cite{Hirsch}). If we denote a
point of $S^p\times S^p$ by $(a,b)$ with $a=(x_1,\ldots,x_{p+1})$
and $b=(x_{p+2},\ldots,x_{2p+2})$, then the restriction of $\Omega$
onto $S^p\times S^p=\partial{\cal H}_i$ will map $(a,b)$ to the
point $({\rm R}(b),a)$ where ${\rm R}(b)$ stands for the ``{\it
reflected"} point $(-x_{p+2},x_{p+3},\ldots,x_{2p+2})$. Such a
diffeomorphism of $S^p\times S^p$ induces in homology the
automorphism of the form $\begin{pmatrix}
0 & -1\\
1 & 0\\
\end{pmatrix}$ and therefore we can define a diffeomorphism $\vf_V$ as
the  restriction $\vf_V:=\Omega|_{S^{2p+2}}$.
\end{proof}

\noindent \underline{Remark}: Notice that $\Omega^4$, as well as
$\Omega^4|_{S^p\times S^p}$ is the identity map. We will use it
later to show that certain short exact sequences split.

\begin{cor}
Let $2\leq p <q$, and $S^p\times S^q$ be standardly embedded in
$S^{p+q+2}$. Then ${\rm Im}(h_{\cal E})\cong\int_2$.
\end{cor}
\begin{proof}
Indeed, it is enough to consider again the union $S^{p+q+1}={\cal
H}_1\cup-{\cal H}_2$ where
$$
{\cal H}_1:=\{(x_1,x_2,\ldots, x_{p+q+2})\in
S^{p+q+1}~|~x^2_1+x^2_2+\cdots+x^2_{p+1}\leq \frac{1}{2}\}
$$ and
${\cal H}_2$ is defined correspondingly and a linear map $\Omega':
\rea^{p+q+3}\lra \rea^{p+q+3}$ defined by the obvious formula
$$
\Omega'(x_0,x_1,\ldots,x_{p+q+1},x_{p+q+2}):=
(x_0,-x_1,x_2\ldots,x_{p+1},-x_{p+2},x_{p+3}\ldots,x_{p+q+2}).
$$
Restriction of this map onto $S^{p+q+2}$ will be an orientation
preserving diffeomorphism that restricts on $S^p\times
S^q=\partial{\cal H}_i$ to a diffeomorphism defined by the formula
$(a,b)\lmt (\rm{R}(a),\rm{R}(b))$, where $a \in S^p$ and $b \in
S^q$. Such a diffeomorphism generates the group ${\rm
Im}(h)\cong\int_2$ (see Proposition 2.1 of \cite{Sato}).
\end{proof}

\begin{cor}
Let $p> 2$ be even, and $S^p\times S^p$ be standardly embedded in
$S^{2p+2}$. Then ${\rm Im}(h_{\cal E})={\rm
Im}(h)\cong\int_2\oplus\int_2$.
\end{cor}
\begin{proof}
This time we consider the orthogonal $(2p+3)\times(2p+3)-$matrix
$$ \widehat{~\Omega} :=
\begin{pmatrix}
\begin{tabular}{c|c|c}
$-1 $ & $\begin{matrix} 0 & 0 & \ldots & 0\end{matrix}$ &
$\begin{matrix} 0 & 0 & \ldots & 0\end{matrix}$ \\
\hline
$\begin{matrix} 0\\
0\\
\vdots \\
0\end{matrix}$ & $
\begin{matrix}
0 & 0 & \ldots & 0 \\
0 & 0 & \ldots & 0 \\
\vdots & \vdots & \ddots & \vdots \\
0 & 0 & \ldots & 0
\end{matrix}$ &
$\begin{matrix}
1 & 0 & \ldots & 0 \\
0 & 1 & \ldots & 0 \\
\vdots & \vdots & \ddots & \vdots \\
0 & 0 & \ldots & 1
\end{matrix}$\\
\hline
$\begin{matrix} 0\\
0\\
\vdots \\
0\end{matrix}$ & $\begin{matrix}
1 & 0 & \ldots & 0 \\
0 & 1 & \ldots & 0 \\
\vdots & \vdots & \ddots & \vdots \\
0 & 0 & \ldots & 1
\end{matrix}$ &
$\begin{matrix}
0 & 0 & \ldots & 0 \\
0 & 0 & \ldots & 0 \\
\vdots & \vdots & \ddots & \vdots \\
0 & 0 & \ldots & 0
\end{matrix}$
\end{tabular}
\end{pmatrix}
$$
which has order two and $\rm{det}(\widehat{~\Omega})=1$. Restriction
of this linear map onto $S^{2p+2}$ will be an orientation preserving
diffeomorphism, which extends the interchange of factors in
$S^p\times S^p$ $\bigl((a,b)\lmt (b,a)$ for $(a,b)\in S^p\times
S^p$, and we assumed again that the equator is $S^{2p+1}={\cal
H}_1\cup-{\cal H}_2$, and $S^p\times S^p=\partial{\cal H}_i\bigr)$.
This interchange generates one of the $\int_2-$copies of ${\rm
Im}(h)$, and the other copy is generated, as above, by the
simultaneous orientation reversal of each factor: $(a,b)\lmt
(\rm{R}(a),\rm{R}(b))$ (cf. \cite{Wall64}, \S1; and \cite{Levine2},
\S1.2).
\end{proof}

\section{Trivial action on the homology}

Consider a ``trivial" $M^n$-knot $(S^{n+2},~M^n)$ with $n\geq 5$
(``triviality" here means that for an equator we have $ S^{n+1} =
W_+\cup W_-$ with $M^n=\partial W_{\pm}$). Identification of the
group $\Theta_{n+1}$ of all homotopy $(n+1)$-spheres with the
relative mapping class group $\pi_0 {\rm Diff}(D^n,{\rm
rel}~\partial)$ (see \cite{Wall25}), gives a homomorphism $\iota:
\Theta_{n+1} \lra \pi_0 {\rm Diff}(M^n)$ (if one presents a homotopy
sphere $\Sigma_{\varphi}$ as the union of two $(n+1)$-disks glued
together via a diffeomorphism $\varphi\in {\rm Diff}(D^n,{\rm
rel}~\partial)$, then $\iota(\Sigma_{\varphi})$ is defined to be the
identity outside an embedded $n$-disk $D^n\hra M^n$ and $\varphi$ on
that disk).

\begin{thm} If a diffeomorphism $\phi\in {\rm Im}(\iota)$
extends to a diffeomorphism $\Phi$ of the trivial $M^n$-knot
$(S^{n+2},~M^n)$, then $\phi$ is pseudo-isotopic to the identity.
\end{thm}
\begin{proof}
Take two canonical pairs of the disks $(D^{n+3}_+,~D^{n+2}_+)$ and
$(D^{n+3}_-,~D^{n+2}_-)$ with the boundary $(S^{n+2},~S^{n+1})$.
Push $W_+$ along the normal vector field inside $D^{n+3}_+$ in such
a way that we get the copy $\widetilde{W}_+$ of $W_+$ embedded in
$D^{n+3}_+$ so that $\widetilde{W}_+ \cap S^{n+2} = M^n$. We can
assume w.l.o.g. that we have $W_+\times I\subset D^{n+3}_+$ with
$W_+\times I\cap S^{n+2} = W_+\times \{0\}$ and $\widetilde{W}_+ =
M\times I\cup W_+\times \{1\}$. Repeat this step with $W_-$ to
obtain $\widetilde{W}_-\subset D^{n+3}_-$. Then we glue together two
pairs $(D^{n+3}_+,~\widetilde{W}_+)$ and
$(D^{n+3}_-,~\widetilde{W}_-)$ via the diffeomorphism $\Phi$ to
obtain a spherical knot
\begin{equation}
\label{glue} (S^{n+3},~\Sigma_{\phi}^{n+1}):=
(D^{n+3}_+,~\widetilde{W}_+) \bigcup_{\Phi}
(D^{n+3}_-,~\widetilde{W}_-).
\end{equation}
We want to show that the homotopy sphere $\Sigma_{\phi}$ bounds a
topological disk inside $S^{n+3}$. Indeed, since $W_+\cup W_- =
S^{n+1}$ bounds a smooth disk $D^{n+2}$ in $S^{n+2}$, it is clear
that the union
$${\cal
D}^{n+2}_1:=W_+\times I \bigcup_{W_+\times \{0\}} D^{n+2} \subset
D^{n+3}_+$$ is homeomorphic to the $(n+2)$-disk with the boundary
$\partial({\cal D}^{n+2}_1) = \widetilde{W}_+\cup_M W_-$. As a
homeomorphism, map $\phi$ is isotopic to the identity, which implies
that the union $\widetilde{W}_-\cup \Phi(W_-)$ (obtained also as a
result of gluing two pairs (\ref{glue}) together) is homeomorphic to
the double ${\cal D}W_-$ (for a manifold $X$ with boundary the
double is defined as the boundary $\partial(X\times I)=:{\cal D}X$).
Hence, we can assume that this union bounds in the ``lower"
hemisphere of $S^{n+3}$ a manifold homeomorphic to the product
$W_-\times I\subset D^{n+3}_-$, and that the intersection of this
product with the equator sphere $S^{n+2}$ is $\Phi(W_-)$. This
implies that the following union
$$
{\cal D}^{n+2}_1 \bigcup_{\Phi(W_-)} W_-\times I =: {\cal D}^{n+2}
$$
will be a topological disk in the sphere $S^{n+3}$ with the boundary
$$
\partial {\cal D}^{n+2}= \widetilde{W}_+\bigcup_{\phi}
\widetilde{W}_- = \Sigma_{\phi}.
$$ Now it follows from Corollary 4.16 of \cite{RandS} that the
pair $(S^{n+3},~\Sigma_{\phi}^{n+1})$ is combinatorially unknotted
and therefore the complement $S^{n+3}\setminus \Sigma_{\phi}^{n+1}$
has the homotopy type of a circle. Using Theorem III from
\cite{Kervaire2} (cf. also \cite{Levine11}) we conclude that the
exotic sphere $\Sigma^{n+1}_{\phi}$ is diffeomorphic to the standard
one. This implies that $\phi$ must be isotopic to the identity (cf.
\cite{Sato},  \S 4).
\end{proof}

\begin{cor}
If $n\geq 5$ and $(S^{n+2},~S^n)$ is the unknot, then ${\cal
E}(S^{n+2},~S^n)\cong Id$.
\end{cor}

\noindent \underline{Remark 1}: If the trivial $M^n$-knot is
actually a spherical knot, i.e. $M^n=S^n$, then it is not hard to
see geometrically that the complement $S^{n+3}\setminus
\Sigma_{\phi}^{n+1}$ to the homotopy sphere $\Sigma_{\phi}^{n+1}$
constructed above is a disk bundle over a circle.

\noindent \underline{Remark 2}: It is interesting to note that a
homotopy $n$-sphere $\Sigma^n$ can be smoothly embedded in $S^{n+2}$
if and only if it bounds a parallelizable manifold (see
\cite{Kervaire2}, \S3).

Let us continue now with the following example. Take the product
$S^{p-2}\times S^{p-1}$ standardly embedded in $S^{p-1}\times
S^{p-1}$ (the $(p-2)$-sphere is the equator of the first
$(p-1)$-sphere), where $p\geq 9$ and $p\equiv 6\pmod{8}$, and again
present the sphere $S^{2p-1}=(D^p\times S^{p-1})\cup (S^{p-1}\times
D^p)$ as the union of two handlebodies. Then the group of the
isotopy classes of diffeomorphisms of $S^{p-2}\times S^{p-1}$ which
act trivially on the homology is isomorphic to
$\pi_{p-1}(SO(p-1))\oplus \Theta_{2p-2}$ (see, for example,
\cite{Turner}, Theorem 3.10). Since $p-1\equiv 5\pmod{8}$, we see
from \cite{Kervaire} that $\pi_{p-1}(SO(p-1))\cong\int_2$ and it
follows also from \cite{Kervaire} that the inclusion
$SO(p-1)\stackrel{\tau}{\hra} SO(p)$ induces the trivial
homomorphism $\tau_*: \pi_{p-1}(SO(p-1))\lra \pi_{p-1}(SO(p))$.
Hence there exists a smooth map $\gamma: D^p\lra SO(p)$ which
extends the composition of the generator $g: S^{p-1}\lra SO(p-1)$ of
$\pi_{p-1}(SO(p-1))$ with $\tau$, i.e. $\gamma|_{S^{p-1}} =
\tau\circ g$. Now we repeat the construction we used to prove Lemma
1. and define a diffeomorphism $\Phi$ of the ambient sphere
$S^{2p-1}$ by the formula:
$$
\Phi(x,y):=\left\{
\begin{array}{lcl}
(\tau(g(y))\circ x,~ y)& {\rm if} & (x,y)\in D^p\times S^{p-1}\\
(\gamma(y)\circ x,~y)& {\rm if} & (x,y)\in S^{p-1}\times D^p\\
\end{array}
\right.
$$
The following lemma immediately follows from Corollary 1. and
Theorem 1.

\begin{lem}
Let $S^{p-2}\times S^{p-1}$ be standardly embedded in $S^{2p-1}$,
with $p\geq 9$ and $p\equiv 6\pmod{8}$. Then the following short
exact sequence splits
$$
0\lra \int_2 \lra {\cal E}(S^{2p-1},~S^{p-2}\times S^{p-1})
\stackrel{h_{\cal E}}{\lra}\int_2 \lra 0.$$
\end{lem}

\noindent Next we generalize this example and study diffeomorphisms
of $S^p\times S^p$ which extend over $S^{2p+2}$ and act trivially on
$H_p(S^p\times S^p;\int)$.

Isotopy classes of the diffeomorphisms which act trivially on
$H_p(S^p\times S^p;\int)$ form a normal subgroup of the mapping
class group and will be denoted by $\pi_0 {\rm SDiff}(S^p\times
S^p)$. Theorem 2 of \cite{Kreck} gives the following description of
such a subgroup for $M=S^p\times S^p,~p\geq 3$:
\begin{equation}
\label{seq} 0\lra \Theta_{2p+1} \stackrel{\iota}{\lra} \pi_0 {\rm
SDiff}(M) \stackrel{\chi}{\lra} {{\rm
Hom}\bigl(H_p(M),~S\pi_p(SO(p))\bigr)} \lra 0.
\end{equation}
Here $S\pi_p(SO(p))$ denotes image of the map $S: \pi_p(SO(p))\to
\pi_p(SO(p+1))$ induced by the inclusion, and the homomorphism
$\chi$ is defined (for $p\geq 4$) as follows (cf. Lemma 1 of
\cite{Kreck}): Take any $f\in {\rm SDiff}(M)$, represent $x\in
H_p(M)$ by a sphere $S^p\hra M$ and use an isotopy to make
$f|_{S^p}=Id$. The stable normal bundle $\nu(S^p)\oplus
\varepsilon^1$ of this sphere in $M$ is trivial and therefore the
differential of $f$ gives an element of $\pi_p(SO(p+1))$. It is
easy to see that this element lies in the image of $S$. If $p=6$
we have $S\pi_p(SO(p))=0$, and for all other $p\geq 3$ the groups
$S\pi_p(SO(p))$ are given in the following table (\cite{Kreck}, p.
644):

\parskip=5mm

\begin{tabular}{|c|c|c|c|c|c|c|c|c|}
\hline $p$ (mod 8) & 0 & 1 & 2 & 3 & 4 & 5 & 6 & 7\\
\hline $S\pi_p(SO(p))$ & ~$\int_2\oplus\int_2$~ & ~$\int_2$~ &
~$\int_2$
~ & ~$\int$ ~ & ~ $\int_2$~ & ~ 0 ~ & ~ $\int_2$ ~ & ~$\int$ \\
\hline
\end{tabular}

\noindent When $p$ is odd, the map $S\pi_p(SO(p)) \hra
\pi_p(SO(p+1))\lra \pi_p(SO(p+2))$ is a monomorphism (see
\cite{Wall4}) and therefore the induced map
$$ {\rm Hom}(H_p(M),~S\pi_p(SO(p))) \lra {{\rm
Hom}(H_p(M),~\pi_p(SO))}$$ is a monomorphism too. Hence we can
describe $\chi([f])[S^p]$ as the class of the stable normal bundle
of $S^p\times S^1$ in the mapping torus $M_f$ of the diffeomorphism
$f$ (cf. Lemma 2 of \cite{Kreck} and recall that $f|_{S^p}=Id$).
Assuming that $f$ extends to a diffeomorphism of $S^{n+2}$, which is
isotopic to the identity, we see that $S^p\times S^1\hra M_f\hra
S^{2p+2}\times S^1$. Since the normal bundle of $M_f$ in
$S^{2p+2}\times S^1$ is trivial, and the stable normal bundle of
$S^p\times S^1$ in $S^{2p+2}\times S^1$ is also trivial, it implies
that the stable normal bundle of $S^p\times S^1$ in $M_f$ will be
trivial as well. Thus the homomorphism $\chi$, restricted onto the
isotopy classes of diffeomorphisms of $S^p\times S^p$, which extend
over $S^{2p+2}$ and act trivially on the homology, is the trivial
map for odd $p\geq 5$. This result together with Theorem 1. and the
exact sequence (\ref{seq}) implies that in such case ${\rm
ker}(h_{\cal E})=\{0\}$. If $p=3$, we have $S\pi_p(SO(p))\simeq
\int$ and we can identify ${\rm Hom}(H_p(M),~S\pi_p(SO(p)))$ with
$H^p(M)$. Then one can use the Pontrjagin class of the mapping torus
$M_f$ instead of $\chi$ (see \cite{Kreck}), to obtain the same
conclusion that the kernel of the map $h_{\cal E}$ is trivial when
$p=3$ as well.

Let us now turn to the case when $p$ is even and at least 4, and
denote the kernel of $h_{\cal E}$ by
$$
{\rm ker}(h_{\cal E}):={\cal SE}(S^{2p+2},~S^p\times S^p).
$$
The following computations are well known (cf. \cite{Wall4}):

$$\begin{tabular}{|c|c|c|c|c|}
\hline $p$ (mod 8) & 0 & 2 & 4 & 6\\
\hline $S\pi_p(SO(p))$ & ~$\int_2\oplus\int_2$~ & ~$\int_2$ ~ & ~
$\int_2$~ & ~ $\int_2$\\
\hline $\pi_p(SO(p+1))$ & ~$\int_2\oplus\int_2$~ & ~$\int_2$ ~ & ~
$\int_2$~ & ~ $\int_2$\\
\hline $\pi_p(SO(p+2))$ & ~$\int_2$~ & ~$0$ ~ & ~
$0$~ & ~ $0$\\
\hline
\end{tabular}
$$

\noindent It follows from these computations and the paragraph
above, that in this case, the image of the restriction of $\chi$
onto the subgroup ${\cal SE}(S^{2p+2},~S^p\times S^p)$ is either
zero or $\int_2\oplus\int_2$.

\begin{lem}
If $p$ is even and $p\geq 4$, then ${\cal SE}(S^{2p+2},~S^p\times
S^p)\cong \int_2\oplus\int_2$ and $ {\cal E}(S^{2p+2},~S^p\times
S^p)\cong {\rm D}_8\oplus \int_2 $, where ${\rm D}_8$ stands for the
dihedral group of order 8.
\end{lem}
\begin{proof}
Indeed, as we just mentioned above, the restriction of $\chi$ takes
values in $ {\rm Hom}(H_p(M),~\int_2)\cong\int_2\oplus\int_2$ and to
prove that this restriction is an epimorphism, it is enough to
construct a diffeomorphism of the pair $(S^{2p+2},~S^p\times S^p)$
for each of the generators of $\int_2\oplus\int_2$. We will use
results from \cite{Sato} to do that. Let $g: S^p\lra SO(p+1)$ be a
smooth map, which represents a nontrivial element of the group
$\pi_p(SO(p+1))$ and such that $\tau\circ g$ is homotopic to the
identity ($\tau$ here is the canonical inclusion $SO(p+1)\hra
SO(p+2)$). Homotopy class of this map $g$ has order two for each
even $p$, as we just saw in the table above. It then follows from
Proposition 3.2 of \cite{Sato} that the diffeomorphisms of
$S^p\times S^p$ defined by the formulas
$$
(x,~y)\stackrel{\delta_1}{\lmt} (x,~g(x)\circ y)~~~~~{\rm
and}~~~~~(x,~y)\stackrel{\delta_2}{\lmt}(g(y)\circ x,~y)
$$ with $(x,y)\in S^p\times S^p$ are representatives of the
generators of the group $\int_2\oplus\int_2$. Since $\tau\circ g\sim
Id$, there exists a smooth map $\gamma: D^{p+1}\lra SO(p+2)$ such
that $\gamma|_{\partial D^{p+1}} = \tau\circ g$ and we define an
orientation preserving diffeomorphism of the sphere $S^{2p+2}$ by
the following formula (in the second case the formula should be
modified in the obvious way):
$$
\Phi(x,y):=\left\{
\begin{array}{lcl}
(x,~\tau(g(x))\circ y)& {\rm if} & (x,y)\in S^p\times D^{p+2}\\
(x,~\gamma(x)\circ y)& {\rm if} & (x,y)\in D^{p+1}\times S^{p+1}\\
\end{array}
\right.
$$ Thus we obtain two required diffeomorphisms of the pair
$(S^{2p+2},~S^p\times S^p)$, and also exactness of the following
sequence:
$$
0 \lra \int_2\oplus\int_2 \lra {\cal E}(S^{2p+2},~S^p\times
S^p)\stackrel{h_{\cal E}}{\lra} \int_2\oplus \int_2 \lra 0
$$

Recall that generators of ${\rm Im}(h_{\cal E})$ are the
simultaneous orientation reversal: $(x,y)\lmt
(\rm{R}(x),\rm{R}(y))$, and the interchange of factors:
$u:~(x,y)\lmt (y,x)$. It is not hard to see that the orientation
reversal acts trivially by conjugation on ${\rm ker}(h_{\cal E})$.
Hence one of the $\int_2$-factors of ${\rm Im}(h_{\cal E})$ splits
off. As for the conjugation by $u$, we have
$$
(x,~y)\stackrel{u}{\lmt} (y,~x)\stackrel{\delta_1}{\lmt}
(y,~g(y)\circ x)\stackrel{u}{\lmt} (g(y)\circ x,~y),
$$
that is $\delta_2=u\circ\delta_1\circ u$. It implies that the factor
group ${\cal E}(S^{2p+2},~S^p\times S^p)/_{\int_2}$ admits the
following presentation: $\langle
a,~b,~u~|~a^2=b^2=u^2=e,~ab=ba,~au=ub \rangle.$ We leave it as an
exercise to show that this presentation gives the dihedral group of
order eight.
\end{proof}

\noindent Let us now summarize what has been proved and state

\begin{thm} Let $S^p\times S^p$ be standardly embedded in
$S^{2p+2}$, and $p\geq 3$. Then
$$
{\cal E}(S^{2p+2},~S^p\times S^p)\cong\left\{
\begin{array}{lcl}
\rm{D}_8\oplus\int_2 & {\rm if} & p~{\rm is~even}\\
\Gv2 & {\rm if} & p~{\rm is~odd}
\end{array}
\right.
$$
\end{thm}

\noindent \underline{Remark 1}: The group of pseudo-isotopy classes
of diffeomorphisms of $S^2\times S^2$ is isomorphic to
$\int_2\oplus\int_2$ (cf. \cite{Wall64}, \S2). Our proof of
Corollary 2 implies that each of the generators has a representative
extendable to a diffeomorphism of $S^6$.

\noindent \underline{Remark 2}: It is interesting to note that the
group $ {\cal E}(S^{2p+2},~S^p\times S^p)$ is isomorphic to $\Gv2 $
for all odd $p\geq 1$ and the ``Hopf invariant one" cases (i.e. when
$p=1,~3~{\rm or}~7$) play no special role here (cf. \cite{Levine2},
\S 1.2, for example).

\end{document}